\documentclass{amsart}
\usepackage{amsfonts}

\newtheorem{theorem}{Theorem}
\theoremstyle{plain}

\newtheorem{remark}{Remark}

\numberwithin{equation}{section}
\input{tcilatex}

\begin{document}
\title[Bessel Type Inequalities]{Bessel Type Inequalities for
Non-Orthonormal Families of Vectors in Inner Product Spaces}
\author{Sever S. Dragomir}
\address{School of Computer Science and Mathematics\\
Victoria University\\
PO Box 14428, Melbourne VIC 8001\\
Australia.}
\email{sever.dragomir@vu.edu.au}
\urladdr{http://rgmia.vu.edu.au/dragomir}
\date{July 21, 2005.}

\begin{abstract}
Some sharp Bessel type inequalities for non-orthonormal families of vectors
in inner product spaces are given. Applications for complex numbers are also
provided.
\end{abstract}

\maketitle

\section{Introduction}

Let $\left( H;\left\langle \cdot ,\cdot \right\rangle \right) $ be an inner
product space over the real or complex number field $\mathbb{K}$. If $%
\left\{ e_{i}\right\} _{i\in \left\{ 1,\dots ,n\right\} }$ are orthonormal
vectors in $H,$ i.e., $\left\langle e_{i},e_{j}\right\rangle =\delta _{ij}$
for all $i,j\in \left\{ 1,\dots ,n\right\} ,$ where $\delta _{ij}$ is the
Kronecker delta, then we have the following inequality:%
\begin{equation}
\sum_{j=1}^{n}\left\vert \left\langle x,e_{j}\right\rangle \right\vert
^{2}\leq \left\Vert x\right\Vert ^{2}\quad \text{for any \ }x\in H,
\label{1.1}
\end{equation}%
which is well known in the literature as Bessel's inequality.

In 1941, R.P. Boas \cite{BO} and in 1944, independently, R. Bellman \cite{BE}
proved the following generalisation of Bessel's inequality (see also \cite[%
p. 392]{MPF}):

\begin{theorem}[Boas-Bellman; 1941, 1944]
\label{t1}If $x,y_{1},\dots ,y_{n}$ are vectors in an inner product space $%
\left( H;\left\langle \cdot ,\cdot \right\rangle \right) ,$ then%
\begin{equation}
\sum_{i=1}^{n}\left\vert \left\langle x,y_{i}\right\rangle \right\vert
^{2}\leq \left\Vert x\right\Vert ^{2}\left[ \max_{1\leq i\leq n}\left\Vert
y_{i}\right\Vert ^{2}+\left( \sum_{1\leq i\neq j\leq n}\left\vert
\left\langle y_{i},y_{j}\right\rangle \right\vert ^{2}\right) ^{\frac{1}{2}}%
\right] .  \label{1.2}
\end{equation}
\end{theorem}

In 1971, E. Bombieri \cite{BOM} (see also \cite[p. 394]{MPF}) gave the
following generalisation of Bessel's inequality:

\begin{theorem}[Bombieri, 1971]
\label{t2}Let $x,y_{1},\dots ,y_{n}$ be vectors in $\left( H;\left\langle
\cdot ,\cdot \right\rangle \right) .$ Then%
\begin{equation}
\sum_{i=1}^{n}\left\vert \left\langle x,y_{i}\right\rangle \right\vert
^{2}\leq \left\Vert x\right\Vert ^{2}\max_{1\leq i\leq n}\left\{
\sum_{j=1}^{n}\left\vert \left\langle y_{i},y_{j}\right\rangle \right\vert
\right\} .  \label{1.3}
\end{equation}
\end{theorem}

Another generalisation of Bessel's inequality was obtained by A. Selberg
(see for example \cite[p. 394]{MPF}):

\begin{theorem}[Selberg]
\label{t3}Let $x,y_{1},\dots ,y_{n}$ be vectors in $H$ with $y_{i}\neq 0$
for $i\in \left\{ 1,\dots ,n\right\} .$ Then%
\begin{equation}
\sum_{i=1}^{n}\frac{\left\vert \left\langle x,y_{i}\right\rangle \right\vert
^{2}}{\sum_{j=1}^{n}\left\vert \left\langle y_{i},y_{j}\right\rangle
\right\vert }\leq \left\Vert x\right\Vert ^{2}.  \label{1.4}
\end{equation}
\end{theorem}

In 2003, the author obtained the following Boas-Bellman type inequality \cite%
{DR1}:

\begin{theorem}[Dragomir, 2003]
\label{t4}For any $x,y_{1},\dots ,y_{n}\in H$ one has%
\begin{equation}
\sum_{i=1}^{n}\left\vert \left\langle x,y_{i}\right\rangle \right\vert
^{2}\leq \left\Vert x\right\Vert ^{2}\left\{ \max_{1\leq i\leq n}\left\Vert
y_{i}\right\Vert ^{2}+\left( n-1\right) \max_{1\leq i\neq j\leq n}\left\vert
\left\langle y_{i},y_{j}\right\rangle \right\vert \right\} .  \label{1.5}
\end{equation}
\end{theorem}

We remark that in all inequalities (\ref{1.2}) -- (\ref{1.5}) the case when $%
\left\{ y_{i}\right\} _{i\in \left\{ 1,\dots ,n\right\} }$ is an orthonormal
family produces the classical Bessel's inequality.

A generalisation of the Bombieri result for the pair $\left( p,q\right) $
with $p>1,$ $\frac{1}{p}+\frac{1}{q}=1$ has been obtained by the author in
2003, \cite{DR2}.

\begin{theorem}[Dragomir, 2003]
\label{t5}For any $x,y_{1},\dots ,y_{n}\in H$ \ with not all $\left\langle
x,y_{i}\right\rangle =0$ for $i\in \left\{ 1,\dots ,n\right\} ,$ one has:%
\begin{equation}
\frac{\left( \sum_{i=1}^{n}\left\vert \left\langle x,y_{i}\right\rangle
\right\vert ^{2}\right) ^{2}}{\left( \sum_{i=1}^{n}\left\vert \left\langle
x,y_{i}\right\rangle \right\vert ^{p}\right) ^{\frac{1}{p}}\left(
\sum_{i=1}^{n}\left\vert \left\langle x,y_{i}\right\rangle \right\vert
^{q}\right) ^{\frac{1}{q}}}\leq \left\Vert x\right\Vert ^{2}\max_{1\leq
i\leq n}\left( \sum_{j=1}^{n}\left\vert \left\langle
y_{i},y_{j}\right\rangle \right\vert \right) .  \label{1.6}
\end{equation}
\end{theorem}

If in this inequality one considers $p=q=2,$ then one obtains Bombieri's
result (\ref{1.3}).

From a different perspective, the following result for the sum of Fourier
coefficients due to H. Heilbronn may be stated \cite{H} (see also \cite[p.
395]{MPF}).

\begin{theorem}[Heilbronn, 1958]
\label{t6}For any $x,y_{1},\dots ,y_{n}\in H$ one has%
\begin{equation}
\sum_{i=1}^{n}\left\vert \left\langle x,y_{i}\right\rangle \right\vert \leq
\left\Vert x\right\Vert \left( \sum_{i,j=1}^{n}\left\vert \left\langle
y_{i},y_{j}\right\rangle \right\vert \right) ^{\frac{1}{2}}.  \label{1.7}
\end{equation}
\end{theorem}

In 1992, Pe\v{c}ari\'{c} \cite{P} (see also \cite[p. 394]{MPF}) proved the
following inequality that incorporates some of the results above:

\begin{theorem}[Pe\v{c}ari\'{c}, 1992]
\label{t7}Let $x,y_{1},\dots ,y_{n}\in H$ and $c_{1},\dots ,c_{n}\in \mathbb{%
K}$. Then%
\begin{align}
\left\vert \sum_{k=1}^{n}c_{k}\left\langle x,y_{k}\right\rangle \right\vert
^{2}& \leq \left\Vert x\right\Vert ^{2}\sum_{i=1}^{n}\left\vert
c_{i}\right\vert ^{2}\left( \sum_{j=1}^{n}\left\vert \left\langle
y_{i},y_{j}\right\rangle \right\vert \right)  \label{1.8} \\
& \leq \left\Vert x\right\Vert ^{2}\sum_{k=1}^{n}\left\vert c_{k}\right\vert
^{2}\max_{1\leq i\leq n}\left\{ \sum_{j=1}^{n}\left\vert \left\langle
y_{i},y_{j}\right\rangle \right\vert \right\} .  \notag
\end{align}
\end{theorem}

He showed that the Bombieri inequality (\ref{1.3}) may be obtained from (\ref%
{1.8}) for the choice $c_{i}=\overline{\left\langle x,y_{i}\right\rangle }$
(using the second inequality), the Selberg inequality (\ref{1.4}) may be
obtained from the first part of (\ref{1.8}) for the choice%
\begin{equation*}
c_{i}=\frac{\overline{\left\langle x,y_{i}\right\rangle }}{%
\sum_{j=1}^{n}\left\vert \left\langle y_{i},y_{j}\right\rangle \right\vert }%
,\quad i\in \left\{ 1,\dots ,n\right\} ,
\end{equation*}%
while the Heilbronn inequality (\ref{1.7}) may be obtained from the first
part of (\ref{1.8}) if one chooses%
\begin{equation*}
c_{i}=\frac{\overline{\left\langle x,y_{i}\right\rangle }}{\left\vert
\left\langle x,y_{i}\right\rangle \right\vert },\quad i\in \left\{ 1,\dots
,n\right\} .
\end{equation*}

In the spirit of Pe\v{c}ari\'{c}'s result, the author proved in \cite{DR2}
the following result as well:

\begin{theorem}[Dragomir, 2004]
\label{t8}Let $x,y_{1},\dots ,y_{n}\in H$ and $c_{1},\dots ,c_{n}\in \mathbb{%
K}$. Then%
\begin{multline}
\left\vert \sum_{k=1}^{n}c_{k}\left\langle x,y_{k}\right\rangle \right\vert
^{2}  \label{1.9} \\
\leq \left\Vert x\right\Vert ^{2}\times \left\{ 
\begin{array}{l}
\max\limits_{1\leq i\leq n}\left\vert c_{k}\right\vert
\sum\limits_{k=1}^{n}\left\vert c_{k}\right\vert \max\limits_{1\leq i\leq
n}\left( \sum\limits_{j=1}^{n}\left\vert \left\langle
y_{i},y_{j}\right\rangle \right\vert \right) ; \\ 
\\ 
\sum\limits_{k=1}^{n}\left\vert c_{k}\right\vert \left(
\sum\limits_{k=1}^{n}\left\vert c_{k}\right\vert ^{p}\right) ^{\frac{1}{p}%
}\max\limits_{1\leq i\leq n}\left\{ \left( \sum\limits_{j=1}^{n}\left\vert
\left\langle y_{i},y_{j}\right\rangle \right\vert ^{q}\right) ^{\frac{1}{q}%
}\right\} ,\  \\ 
\hfill p>1,\ \frac{1}{p}+\frac{1}{q}=1; \\ 
\\ 
\left( \sum\limits_{k=1}^{n}\left\vert c_{k}\right\vert \right)
^{2}\max\limits_{1\leq i,j\leq n}\left\vert \left\langle
y_{i},y_{j}\right\rangle \right\vert .%
\end{array}%
\right.
\end{multline}
\end{theorem}

In particular, for the choice $c_{k}=\overline{\left\langle
x,y_{k}\right\rangle },$ $k\in \left\{ 1,\dots ,n\right\} ,$ we deduce from (%
\ref{1.9}) that%
\begin{equation}
\frac{\left( \sum_{k=1}^{n}\left\vert \left\langle x,y_{k}\right\rangle
\right\vert ^{2}\right) ^{2}}{\max\limits_{1\leq k\leq n}\left\vert
\left\langle x,y_{k}\right\rangle \right\vert \sum_{k=1}^{n}\left\vert
\left\langle x,y_{k}\right\rangle \right\vert }\leq \left\Vert x\right\Vert
^{2}\max\limits_{1\leq i\leq n}\left\{ \sum\limits_{j=1}^{n}\left\vert
\left\langle y_{i},y_{j}\right\rangle \right\vert \right\} ,  \label{1.10}
\end{equation}%
\begin{gather}
\frac{\left( \sum_{k=1}^{n}\left\vert \left\langle x,y_{k}\right\rangle
\right\vert ^{2}\right) ^{2}}{\sum_{k=1}^{n}\left\vert \left\langle
x,y_{k}\right\rangle \right\vert \left( \sum_{k=1}^{n}\left\vert
\left\langle x,y_{k}\right\rangle \right\vert ^{p}\right) ^{\frac{1}{p}}}%
\leq \left\Vert x\right\Vert ^{2}\max\limits_{1\leq i\leq n}\left\{ \left(
\sum\limits_{j=1}^{n}\left\vert \left\langle y_{i},y_{j}\right\rangle
\right\vert ^{q}\right) ^{\frac{1}{q}}\right\} ,  \label{1.11} \\
\notag
\end{gather}%
for \ $p>1,\ \frac{1}{p}+\frac{1}{q}=1;$ and%
\begin{equation}
\frac{\left( \sum_{k=1}^{n}\left\vert \left\langle x,y_{k}\right\rangle
\right\vert ^{2}\right) ^{2}}{\left( \sum_{k=1}^{n}\left\vert \left\langle
x,y_{k}\right\rangle \right\vert \right) ^{2}}\leq \left\Vert x\right\Vert
^{2}\max\limits_{1\leq i,j\leq n}\left\vert \left\langle
y_{i},y_{j}\right\rangle \right\vert ,  \label{1.12}
\end{equation}%
provided not all $\left\langle x,y_{k}\right\rangle ,$ $k\in \left\{ 1,\dots
,n\right\} $ are zero.

The aim of the present paper is to provide different upper bounds for the
Bessel sum $\sum_{i=1}^{n}\left\vert \left\langle x,y_{i}\right\rangle
\right\vert ^{2}$ under various conditions for the vectors enclosed.
Applications for complex numbers are provided as well.

\section{The Results}

The following sharp inequality of Bessel type for Fourier coefficients
satisfying some restrictions may be stated:

\begin{theorem}
\label{t2.1}Let $\left( H;\left\langle \cdot ,\cdot \right\rangle \right) $
be a real or complex inner product space, $\gamma ,\Gamma \in \mathbb{K}$
with $\Gamma \neq -\gamma $ and $x,y_{j}\in H,$ $j\in \left\{ 1,\dots
,n\right\} $ such that%
\begin{multline}
\left( \func{Re}\Gamma -\func{Re}\left\langle x,y_{j}\right\rangle \right)
\left( \func{Re}\left\langle x,y_{j}\right\rangle -\func{Re}\gamma \right) 
\label{1} \\
+\left( \func{Im}\Gamma -\func{Im}\left\langle x,y_{j}\right\rangle \right)
\left( \func{Im}\left\langle x,y_{j}\right\rangle -\func{Im}\gamma \right)
\geq 0
\end{multline}%
or, equivalently,%
\begin{equation}
\left\vert \left\langle x,y_{j}\right\rangle -\frac{\gamma +\Gamma }{2}%
\right\vert \leq \frac{1}{2}\left\vert \Gamma -\gamma \right\vert   \label{2}
\end{equation}%
for each \ $j\in \left\{ 1,\dots ,n\right\} .$ Then%
\begin{equation}
\left( \sum\limits_{j=1}^{n}\left\vert \left\langle x,y_{j}\right\rangle
\right\vert ^{2}\right) ^{\frac{1}{2}}\leq \frac{1}{\sqrt{n}}\left\Vert
x\right\Vert \left\Vert \sum\limits_{j=1}^{n}y_{j}\right\Vert +\frac{1}{4}%
\sqrt{n}\cdot \frac{\left\vert \Gamma -\gamma \right\vert ^{2}}{\left\vert
\Gamma +\gamma \right\vert }.  \label{3}
\end{equation}%
The equality holds in (\ref{3}) for $x\neq 0$ if and only if the equality
case holds in (\ref{2}) for each $j\in \left\{ 1,\dots ,n\right\} $ and%
\begin{equation}
\frac{1}{n}\sum\limits_{j=1}^{n}y_{j}=\frac{1}{4}\cdot \frac{\left\vert
\Gamma \right\vert ^{2}+6\func{Re}\left( \Gamma \bar{\gamma}\right)
+\left\vert \gamma \right\vert ^{2}}{\left( \Gamma +\gamma \right)
\left\Vert x\right\Vert ^{2}}\cdot x.  \label{4}
\end{equation}
\end{theorem}

\begin{proof}
The equivalence between (\ref{1}) and (\ref{2}) is obvious since for the
complex numbers $z,\gamma ,\Gamma $ the following statements are equivalent:

\begin{enumerate}
\item[(i)] $\func{Re}\left[ \left( \Gamma -z\right) \left( \bar{z}-\bar{%
\gamma}\right) \right] \geq 0;$

\item[(ii)] $\left\vert z-\frac{\gamma +\Gamma }{2}\right\vert \leq \frac{1}{%
2}\left\vert \Gamma -\gamma \right\vert .$
\end{enumerate}

The inequality (\ref{2}) is clearly equivalent to%
\begin{equation}
\left\vert \left\langle x,y_{j}\right\rangle \right\vert ^{2}+\left\vert 
\frac{\gamma +\Gamma }{2}\right\vert ^{2}\leq \frac{1}{4}\left\vert \Gamma
-\gamma \right\vert ^{2}+\func{Re}\left[ \left( \bar{\Gamma}+\bar{\gamma}%
\right) \left\langle x,y_{j}\right\rangle \right]  \label{5}
\end{equation}%
for $j\in \left\{ 1,\dots ,n\right\} $ with equality iff the case of
equality is realised in (\ref{2}).

Summing over $j$ from 1 to $n$ in (\ref{5}), we get:%
\begin{equation}
\sum\limits_{j=1}^{n}\left\vert \left\langle x,y_{j}\right\rangle
\right\vert ^{2}+n\left\vert \frac{\gamma +\Gamma }{2}\right\vert ^{2}\leq 
\frac{1}{4}n\left\vert \Gamma -\gamma \right\vert ^{2}+\func{Re}\left[
\left( \bar{\Gamma}+\bar{\gamma}\right) \left\langle
x,\sum\limits_{j=1}^{n}y_{j}\right\rangle \right]  \label{6}
\end{equation}%
with equality if and only if the equality case holds for each $j$ in (\ref{2}%
).

Utilising the arithmetic mean -- geometric mean inequality we can state%
\begin{equation}
2\sqrt{n}\left\vert \frac{\gamma +\Gamma }{2}\right\vert \left(
\sum\limits_{j=1}^{n}\left\vert \left\langle x,y_{j}\right\rangle
\right\vert ^{2}\right) ^{\frac{1}{2}}\leq \sum\limits_{j=1}^{n}\left\vert
\left\langle x,y_{j}\right\rangle \right\vert ^{2}+n\left\vert \frac{\gamma
+\Gamma }{2}\right\vert ^{2}  \label{7}
\end{equation}%
with equality if and only if%
\begin{equation}
\sum\limits_{j=1}^{n}\left\vert \left\langle x,y_{j}\right\rangle
\right\vert ^{2}=n\cdot \left\vert \frac{\gamma +\Gamma }{2}\right\vert ^{2}.
\label{8}
\end{equation}%
Combining (\ref{6}) with (\ref{7}) we deduce%
\begin{align}
\left( \sum\limits_{j=1}^{n}\left\vert \left\langle x,y_{j}\right\rangle
\right\vert ^{2}\right) ^{\frac{1}{2}}& \leq \frac{1}{4}\sqrt{n}\cdot \frac{%
\left\vert \Gamma -\gamma \right\vert ^{2}}{\left\vert \Gamma +\gamma
\right\vert }+\func{Re}\left[ \frac{\left( \bar{\Gamma}+\bar{\gamma}\right) 
}{\Gamma +\gamma }\left\langle x,\sum\limits_{j=1}^{n}y_{j}\right\rangle %
\right]  \label{9} \\
& \leq \frac{1}{4}\sqrt{n}\cdot \frac{\left\vert \Gamma -\gamma \right\vert
^{2}}{\left\vert \Gamma +\gamma \right\vert }+\left\vert \left\langle
x,\sum\limits_{j=1}^{n}y_{j}\right\rangle \right\vert  \notag \\
& \leq \frac{1}{4}\sqrt{n}\cdot \frac{\left\vert \Gamma -\gamma \right\vert
^{2}}{\left\vert \Gamma +\gamma \right\vert }+\left\Vert x\right\Vert
\left\Vert \sum\limits_{j=1}^{n}y_{j}\right\Vert .  \notag
\end{align}%
For the last inequality in (\ref{9}) we have used Schwarz's inequality $%
\left\vert \left\langle u,v\right\rangle \right\vert \leq \left\Vert
u\right\Vert \left\Vert v\right\Vert ,$ $u,v\in H,$ for which, since%
\begin{equation*}
\left\Vert u-\frac{\left\langle u,v\right\rangle v}{\left\Vert v\right\Vert
^{2}}\right\Vert ^{2}=\frac{\left\Vert u\right\Vert ^{2}\left\Vert
v\right\Vert ^{2}-\left\vert \left\langle u,v\right\rangle \right\vert ^{2}}{%
\left\Vert v\right\Vert ^{2}},\quad v\neq 0,
\end{equation*}%
the equality case holds if and only if%
\begin{equation*}
u=\frac{\left\langle u,v\right\rangle v}{\left\Vert v\right\Vert ^{2}}.
\end{equation*}%
Therefore the equality case holds in the last part of (\ref{9}) if and only
if%
\begin{equation*}
\sum\limits_{j=1}^{n}y_{j}=\frac{\sum_{j=1}^{n}\overline{\left\langle
x,y_{j}\right\rangle }x}{\left\Vert x\right\Vert ^{2}}.
\end{equation*}%
Now, if the equality case holds in (\ref{2}) for each $j\in \left\{ 1,\dots
,n\right\} ,$ then squaring and summing over $j\in \left\{ 1,\dots
,n\right\} ,$ we deduce:%
\begin{equation}
\sum\limits_{j=1}^{n}\left\vert \left\langle x,y_{j}\right\rangle
\right\vert ^{2}=\func{Re}\left\langle x,\left( \gamma +\Gamma \right)
\sum\limits_{j=1}^{n}y_{j}\right\rangle +\frac{1}{4}n\left\vert \Gamma
-\gamma \right\vert ^{2}-\frac{1}{4}n\left\vert \Gamma +\gamma \right\vert
^{2}.  \label{10}
\end{equation}%
From (\ref{4}), taking the inner product and the real part we have%
\begin{align}
\func{Re}\left\langle x,\left( \gamma +\Gamma \right)
\sum\limits_{j=1}^{n}y_{j}\right\rangle & =\frac{n}{4}\cdot \frac{\left\vert
\Gamma \right\vert ^{2}+6\func{Re}\left( \Gamma \bar{\gamma}\right)
+\left\vert \gamma \right\vert ^{2}}{\left\Vert x\right\Vert ^{2}}\left\Vert
x\right\Vert ^{2}  \label{11} \\
& =\frac{n}{4}\left[ 2\left\vert \Gamma +\gamma \right\vert ^{2}-\left\vert
\Gamma -\gamma \right\vert ^{2}\right] .  \notag
\end{align}%
Therefore, by (\ref{10}) and (\ref{11}) we get%
\begin{align}
\sum\limits_{j=1}^{n}\left\vert \left\langle x,y_{j}\right\rangle
\right\vert ^{2}& =\frac{1}{4}n\left[ 2\left\vert \Gamma +\gamma \right\vert
^{2}-\left\vert \Gamma -\gamma \right\vert ^{2}\right] +\frac{1}{4}%
n\left\vert \Gamma -\gamma \right\vert ^{2}-\frac{1}{4}n\left\vert \Gamma
+\gamma \right\vert ^{2}  \label{12} \\
& =\frac{n}{4}\left\vert \Gamma +\gamma \right\vert ^{2}.  \notag
\end{align}%
Taking the norm in (\ref{4}) we have%
\begin{align}
& \frac{1}{\sqrt{n}}\left\Vert x\right\Vert \left\Vert
\sum\limits_{j=1}^{n}y_{j}\right\Vert +\frac{1}{4}\sqrt{n}\cdot \frac{%
\left\vert \Gamma -\gamma \right\vert ^{2}}{\left\vert \Gamma +\gamma
\right\vert }  \label{13} \\
& =\frac{\sqrt{n}}{4}\cdot \frac{\left( 2\left\vert \Gamma +\gamma
\right\vert ^{2}-\left\vert \Gamma -\gamma \right\vert \right) }{\left\vert
\Gamma +\gamma \right\vert }+\frac{1}{4}\sqrt{n}\cdot \frac{\left\vert
\Gamma -\gamma \right\vert ^{2}}{\left\vert \Gamma +\gamma \right\vert } 
\notag \\
& =\frac{\sqrt{n}}{2}\left\vert \Gamma +\gamma \right\vert .  \notag
\end{align}%
The equations (\ref{12}) and (\ref{13}) show that the equality case is
realised in (\ref{3}). Conversely, if the equality case holds in (\ref{3}),
it must hold in all inequalities needed to prove it, therefore, we must have:%
\begin{gather}
\left\vert \left\langle x,y_{j}\right\rangle -\frac{\Gamma +\gamma }{2}%
\right\vert =\frac{1}{2}\left\vert \Gamma -\gamma \right\vert \quad \text{%
for each \ }j\in \left\{ 1,\dots ,n\right\} ,  \label{14} \\
\sum\limits_{j=1}^{n}\left\vert \left\langle x,y_{j}\right\rangle
\right\vert ^{2}=n\left\vert \frac{\gamma +\Gamma }{2}\right\vert ^{2},
\label{15} \\
\func{Im}\left\langle x,\left( \Gamma +\gamma \right)
\sum\limits_{j=1}^{n}y_{j}\right\rangle =0  \label{16}
\end{gather}%
and%
\begin{equation}
\sum\limits_{j=1}^{n}y_{j}=\frac{\overline{\left\langle
x,\sum_{j=1}^{n}y_{j}\right\rangle }}{\left\Vert x\right\Vert ^{2}}\cdot x.
\label{17}
\end{equation}%
From (\ref{14}) we get:%
\begin{equation*}
\func{Re}\left\langle x,\left( \Gamma +\gamma \right) y_{j}\right\rangle
=\left\vert \left\langle x,y_{j}\right\rangle \right\vert ^{2}+\frac{1}{4}%
\left\vert \Gamma +\gamma \right\vert ^{2}-\frac{1}{4}\left\vert \Gamma
-\gamma \right\vert ^{2},
\end{equation*}%
which, by summation over $j$ and (\ref{15}), gives%
\begin{align}
\func{Re}\left\langle x,\left( \Gamma +\gamma \right)
\sum\limits_{j=1}^{n}y_{j}\right\rangle & =n\left[ \left\vert \frac{\Gamma
+\gamma }{2}\right\vert ^{2}-\left\vert \frac{\Gamma -\gamma }{2}\right\vert
^{2}\right]  \label{18} \\
& =\frac{n}{4}\cdot \left[ \left\vert \Gamma \right\vert ^{2}+6\func{Re}%
\left( \Gamma \bar{\gamma}\right) +\left\vert \gamma \right\vert ^{2}\right]
.  \notag
\end{align}%
On multiplying (\ref{17}) by $\gamma +\Gamma \neq 0,$ we have%
\begin{align}
\left( \Gamma +\gamma \right) \cdot \sum\limits_{j=1}^{n}y_{j}& =\frac{%
\overline{\left\langle x,\left( \Gamma +\gamma \right)
\sum_{j=1}^{n}y_{j}\right\rangle }}{\left\Vert x\right\Vert ^{2}}\cdot x
\label{19} \\
& =\frac{\func{Re}\left\langle x,\left( \Gamma +\gamma \right)
\sum_{j=1}^{n}y_{j}\right\rangle -i\func{Im}\left\langle x,\left( \Gamma
+\gamma \right) \sum_{j=1}^{n}y_{j}\right\rangle }{\left\Vert x\right\Vert
^{2}}\cdot x.  \notag
\end{align}%
Finally, on making use of (\ref{16}), (\ref{18}) and (\ref{19}), we deduce
the equality (\ref{4}) and the proof is complete.
\end{proof}

The following results that provide a different bound for the Bessel sum $%
\sum_{j=1}^{n}\left\vert \left\langle x,y_{j}\right\rangle \right\vert ^{2}$
may be stated as well.

\begin{theorem}
\label{t2.2}Let $\left( H;\left\langle \cdot ,\cdot \right\rangle \right) $
be a real or complex inner product space, $\Gamma ,\gamma \in \mathbb{K}$
with $\func{Re}\left( \Gamma \bar{\gamma}\right) >0$ and $x,y_{j}\in H,$ $%
j\in \left\{ 1,\dots ,n\right\} $ such that either (\ref{1}) or (\ref{2})
hold true. Then%
\begin{equation}
\sum\limits_{j=1}^{n}\left\vert \left\langle x,y_{j}\right\rangle
\right\vert ^{2}\leq \frac{1}{n}\cdot \frac{\left\vert \Gamma +\gamma
\right\vert ^{2}}{4\func{Re}\left( \Gamma \bar{\gamma}\right) }\left\Vert
\sum\limits_{j=1}^{n}y_{j}\right\Vert ^{2}\left\Vert x\right\Vert ^{2}.
\label{20}
\end{equation}%
The equality holds in (\ref{20}) for $x\neq 0$ if and only if the equality
case holds in (\ref{2}) for each $j\in \left\{ 1,\dots ,n\right\} $ and%
\begin{equation}
\frac{1}{n}\cdot \sum\limits_{j=1}^{n}y_{j}=\frac{2\func{Re}\left( \Gamma 
\bar{\gamma}\right) }{\left( \Gamma +\gamma \right) \left\Vert x\right\Vert
^{2}}\cdot x.  \label{21}
\end{equation}
\end{theorem}

\begin{proof}
From (\ref{6}) we have%
\begin{equation}
\sum\limits_{j=1}^{n}\left\vert \left\langle x,y_{j}\right\rangle
\right\vert ^{2}+n\func{Re}\left( \Gamma \bar{\gamma}\right) \leq \func{Re}%
\left[ \left( \bar{\Gamma}+\bar{\gamma}\right) \left\langle
x,\sum\limits_{j=1}^{n}y_{j}\right\rangle \right]  \label{21.a}
\end{equation}%
with equality if and only if the case of equality holds in (\ref{2}) for
each $j\in \left\{ 1,\dots ,n\right\} .$

Utilising the elementary inequality between the arithmetic and geometric
mean, we have%
\begin{equation}
2\sqrt{n}\left( \sum\limits_{j=1}^{n}\left\vert \left\langle
x,y_{j}\right\rangle \right\vert ^{2}\right) ^{\frac{1}{2}}\sqrt{\func{Re}%
\left( \Gamma \bar{\gamma}\right) }\leq \sum\limits_{i=1}^{n}\left\vert
\left\langle x,e_{i}\right\rangle \right\vert ^{2}+n\func{Re}\left( \Gamma 
\bar{\gamma}\right)  \label{22}
\end{equation}%
with equality if and only if%
\begin{equation}
\sum\limits_{j=1}^{n}\left\vert \left\langle x,y_{j}\right\rangle
\right\vert ^{2}=n\func{Re}\left( \Gamma \bar{\gamma}\right) .  \label{23}
\end{equation}%
Combining (\ref{21.a}) with (\ref{22}) we deduce%
\begin{align}
\sum\limits_{j=1}^{n}\left\vert \left\langle x,y_{j}\right\rangle
\right\vert ^{2}& \leq \frac{\left\{ \func{Re}\left[ \left( \bar{\Gamma}+%
\bar{\gamma}\right) \left\langle x,\sum_{j=1}^{n}y_{j}\right\rangle \right]
\right\} ^{2}}{4n\func{Re}\left( \Gamma \bar{\gamma}\right) }  \label{24} \\
& \leq \frac{\left\vert \Gamma +\gamma \right\vert ^{2}\left\vert
\left\langle x,\sum_{j=1}^{n}y_{j}\right\rangle \right\vert ^{2}}{4n\func{Re}%
\left( \Gamma \bar{\gamma}\right) }  \notag \\
& \leq \frac{1}{n}\cdot \frac{\left\vert \Gamma +\gamma \right\vert ^{2}}{4%
\func{Re}\left( \Gamma \bar{\gamma}\right) }\left\Vert
\sum\limits_{j=1}^{n}y_{j}\right\Vert ^{2}\left\Vert x\right\Vert ^{2}, 
\notag
\end{align}%
where, for the last inequality we have used Schwarz's inequality. The
equality case holds in the last inequality (\ref{24}) iff%
\begin{equation}
\sum\limits_{j=1}^{n}y_{j}=\frac{\overline{\left\langle
x,\sum_{j=1}^{n}y_{j}\right\rangle }\ x}{\left\Vert x\right\Vert ^{2}}.
\label{25}
\end{equation}%
Now, if the equality case holds in (\ref{2}) for each $j\in \left\{ 1,\dots
,n\right\} ,$ then squaring and summing over $j\in \left\{ 1,\dots
,n\right\} $ we deduce:%
\begin{equation}
\sum\limits_{j=1}^{n}\left\vert \left\langle x,y_{j}\right\rangle
\right\vert ^{2}=\func{Re}\left\langle x,\left( \Gamma +\gamma \right)
\sum_{j=1}^{n}y_{j}\right\rangle +\frac{1}{4}n\left\vert \Gamma -\gamma
\right\vert ^{2}-\frac{1}{4}n\left\vert \Gamma +\gamma \right\vert ^{2}.
\label{26}
\end{equation}%
From (\ref{21}), taking the inner product and the real part we have%
\begin{equation}
\func{Re}\left\langle x,\left( \Gamma +\gamma \right)
\sum_{j=1}^{n}y_{j}\right\rangle =2n\func{Re}\left( \Gamma \bar{\gamma}%
\right) .  \label{27}
\end{equation}%
Therefore, by (\ref{26}) and (\ref{27}) we have%
\begin{equation*}
\sum\limits_{j=1}^{n}\left\vert \left\langle x,y_{j}\right\rangle
\right\vert ^{2}=n\func{Re}\left( \Gamma \bar{\gamma}\right) .
\end{equation*}%
Taking the norm in (\ref{21}) we have%
\begin{equation*}
\frac{1}{n}\cdot \frac{\left\vert \Gamma +\gamma \right\vert ^{2}}{4\func{Re}%
\left( \Gamma \bar{\gamma}\right) }\left\Vert
\sum\limits_{j=1}^{n}y_{j}\right\Vert ^{2}\left\Vert x\right\Vert ^{2}=n%
\func{Re}\left( \Gamma \bar{\gamma}\right) .
\end{equation*}%
showing that the equality case holds true in (\ref{20}).

Conversely, if the equality case holds in (\ref{20}), it must hold in all
inequalities needed to prove it, therefore we must have (\ref{14}), (\ref{23}%
), (\ref{25}) and%
\begin{equation}
\func{Im}\left\langle x,\left( \Gamma +\gamma \right)
\sum_{j=1}^{n}y_{j}\right\rangle =0  \label{27.a}
\end{equation}

From (\ref{14}) we have%
\begin{equation*}
\func{Re}\left\langle x,\left( \Gamma +\gamma \right) y_{j}\right\rangle
=\left\vert \left\langle x,y_{j}\right\rangle \right\vert ^{2}+\frac{1}{4}%
\left\vert \Gamma +\gamma \right\vert ^{2}-\frac{1}{4}\left\vert \Gamma
-\gamma \right\vert ^{2},
\end{equation*}%
which by summation over $j$ and (\ref{23}) gives%
\begin{equation}
\func{Re}\left\langle x,\left( \Gamma +\gamma \right)
\sum_{j=1}^{n}y_{j}\right\rangle =n\func{Re}\left( \Gamma \bar{\gamma}%
\right) +n\func{Re}\left( \Gamma \bar{\gamma}\right) =2n\func{Re}\left(
\Gamma \bar{\gamma}\right) .  \label{27.b}
\end{equation}%
Since $\gamma +\Gamma $ is not zero (because $\left\vert \Gamma +\gamma
\right\vert ^{2}\geq 4\func{Re}\left( \Gamma \overline{\gamma }\right) >0),$
we have by (\ref{25}), (\ref{27.a}) and (\ref{27.b}) that%
\begin{align}
\sum_{j=1}^{n}y_{j}& =\frac{\func{Re}\left\langle x,\left( \Gamma +\gamma
\right) \sum_{j=1}^{n}y_{j}\right\rangle -i\func{Im}\left\langle x,\left(
\Gamma +\gamma \right) \sum_{j=1}^{n}y_{j}\right\rangle }{\left( \Gamma
+\gamma \right) \left\Vert x\right\Vert ^{2}}\cdot x  \label{28} \\
& =\frac{2n\func{Re}\left( \Gamma \bar{\gamma}\right) }{\left( \Gamma
+\gamma \right) \left\Vert x\right\Vert ^{2}}\cdot x,  \notag
\end{align}%
which obviously imply the identity (\ref{21}), and the proof of the theorem
is complete.
\end{proof}

\begin{remark}
A more convenient sufficient condition for (\ref{1}) to hold is%
\begin{equation}
\func{Re}\Gamma \geq \func{Re}\left\langle x,y_{j}\right\rangle \geq \func{Re%
}\gamma \quad \text{and}\quad \func{Im}\Gamma \geq \func{Im}\left\langle
x,y_{j}\right\rangle \geq \func{Im}\gamma  \label{29}
\end{equation}%
for each $j\in \left\{ 1,\dots ,n\right\} .$
\end{remark}

\begin{remark}
If $\left\{ e_{1},\dots ,e_{n}\right\} $ is an orthonormal family of
vectors, then from (\ref{3}) we get%
\begin{equation}
\left( \sum\limits_{j=1}^{n}\left\vert \left\langle x,e_{j}\right\rangle
\right\vert ^{2}\right) ^{\frac{1}{2}}\leq \left\Vert x\right\Vert +\frac{1}{%
4}\sqrt{n}\frac{\left\vert \Gamma -\gamma \right\vert ^{2}}{\left\vert
\Gamma +\gamma \right\vert }  \label{30}
\end{equation}%
while from (\ref{20}) we get%
\begin{equation}
\sum\limits_{j=1}^{n}\left\vert \left\langle x,e_{j}\right\rangle
\right\vert ^{2}\leq \frac{\left\vert \Gamma +\gamma \right\vert }{4\func{Re}%
\left( \Gamma \bar{\gamma}\right) }\left\Vert x\right\Vert ^{2}.  \label{31}
\end{equation}%
One must observe that in this case both (\ref{30}) and (\ref{31}) provide
coarser bounds than Bessel's inequality (\ref{1.1}). Therefore, since (\ref%
{3}) and (\ref{20}) are sharp inequalities, they must be used only in the
case of non-orthogonal vectors satisfying (\ref{1}) or (\ref{2}).
\end{remark}

\section{Applications for Complex Numbers}

Utilising Theorems \ref{t2.1} and \ref{t2.2} above, one can sate the
following reverses of the generalised triangle inequality for complex
numbers that may be of interest in applications:

If 
\begin{equation}
\left( \func{Re}\Gamma -\func{Re}z_{j}\right) \left( \func{Re}z_{j}-\func{Re}%
\gamma \right) +\left( \func{Im}\Gamma -\func{Im}z_{j}\right) \left( \func{Im%
}z_{j}-\func{Im}\gamma \right) \geq 0  \label{3.1.1}
\end{equation}%
or, equivalently,%
\begin{equation}
\left\vert z_{j}-\frac{\gamma +\Gamma }{2}\right\vert \leq \frac{1}{2}%
\left\vert \Gamma -\gamma \right\vert   \label{3.2}
\end{equation}%
for each \ $j\in \left\{ 1,\dots ,n\right\} ,$ then

\begin{equation}
\left( \sum\limits_{j=1}^{n}\left\vert z_{j}\right\vert ^{2}\right)
^{1/2}\leq \frac{1}{\sqrt{n}}\left\vert
\sum\limits_{j=1}^{n}z_{j}\right\vert +\frac{1}{4}\sqrt{n}\frac{\left\vert
\Gamma -\gamma \right\vert ^{2}}{\left\vert \Gamma +\gamma \right\vert }
\label{3.3}
\end{equation}%
provided $\Gamma \neq -\gamma ,$ and%
\begin{equation}
\sum\limits_{j=1}^{n}\left\vert z_{j}\right\vert ^{2}\leq \frac{1}{n}\cdot 
\frac{\left\vert \Gamma +\gamma \right\vert ^{2}}{4\func{Re}\left( \Gamma 
\bar{\gamma}\right) }\left\vert \sum\limits_{j=1}^{n}z_{j}\right\vert ^{2}
\label{3.4}
\end{equation}%
provided $\func{Re}\left( \Gamma \bar{\gamma}\right) >0.$

The equality holds in (\ref{3.3}) if and only if the equality case holds in (%
\ref{3.1.1}) (or in (\ref{3.2})) for each $j\in \left\{ 1,\dots ,n\right\} $
and%
\begin{equation*}
\frac{1}{n}\sum\limits_{j=1}^{n}z_{j}=\frac{1}{4}\cdot \frac{\left\vert
\Gamma \right\vert ^{2}+6\func{Re}\left( \Gamma \bar{\gamma}\right)
+\left\vert \gamma \right\vert ^{2}}{\overline{\Gamma }+\overline{\gamma }}.
\end{equation*}%
The equality holds in (\ref{3.4}) if and only if the equality case holds in (%
\ref{3.1.1}) (or in (\ref{3.2})) for each $j\in \left\{ 1,\dots ,n\right\} $
and%
\begin{equation*}
\frac{1}{n}\sum\limits_{j=1}^{n}z_{j}=\frac{2\func{Re}\left( \Gamma \bar{%
\gamma}\right) }{\overline{\Gamma }+\overline{\gamma }}.
\end{equation*}

\end{document}